\documentclass[reqno]{amsart}
\usepackage{amsmath,amssymb,amsfonts,hhline}
\usepackage[arrow,matrix,curve]{xy}
\usepackage{latexsym}
\usepackage{amsthm}
\usepackage{amssymb}
\usepackage{amsmath}
\usepackage{amscd}
\usepackage{amsopn}
\usepackage{pictexwd,dcpic}
\usepackage{hyperref}
\usepackage{comment}

\theoremstyle{definition}
\newtheorem{theorem}{Theorem}[section]
\newtheorem{definition}{Definition}[section]

\newtheorem{lemma}{Lemma}

\begin{document}

\title{The Puiseux characteristic of a Goursat germ}
\author{Corey Shanbrom}
\address{Department of Mathematics, UC Santa Cruz,
Santa Cruz, CA 95064, USA}
\email{cshanbro@ucsc.edu}

\date{6 March, 2013}

\begin{abstract}
Germs of Goursat distributions can be classified according to a geometric coding called an RVT code.  Jean (\cite{J}) and Mormul (\cite{Mo1}) have shown that this coding carries precisely the same data as the small growth vector.  Montgomery and Zhitomirskii (\cite{MZ1}) have shown that such germs correspond to finite jets of Legendrian curve germs, and that the RVT coding corresponds to the classical invariant in the singularity theory of planar curves: the Puiseux characteristic.  Here we derive a simple formula, Theorem \ref{maintheorem}, for the Puiseux characteristic of the curve corresponding to a Goursat germ with given small growth vector.  The simplicity of our theorem (compared with the more complex algorithms previously known) suggests a deeper connection between singularity theory and the theory of nonholonomic distributions.
\end{abstract}
\maketitle


\section{Introduction}\label{intro}
In this paper we give a formula for the Puiseux characteristic of an analytic plane curve germ which represents a Goursat distribution germ with prescribed small growth vector.  It has been shown earlier that the Puiseux characteristic and the small growth vector are equivalent data, as both are equivalent to a geometric stratification of Goursat germs known as RVT coding.  (Puiseux characteristics do not cover the entire RVT stratification, but the difference is secondary.)  However, the existing algorithms for calculating one from another are cumbersome and recursive.  
  Here, we effectively compose two known algorithms and derive a tidy formula which greatly simplifies the existing methods. 
In addition, the Puiseux characteristic is a much more compact labelling than the small growth vector or the RVT code, and our formula allows for a quick conversion to this more convenient labelling.
   The problem solved herein was first proposed in \cite{MZ1} as Question 9.19, part 3, and was asked again in \cite{Mo2} in the Afterword.

    Goursat distributions are located in the antipodes of integrable distributions, as they are \emph{bracket-generating}.   Cartan (\cite{C}) studied the model of the canonical contact distribution on the jet space $J^{k}(\mathbb{R}, \mathbb{R})$. All Goursat distributions were believed to be equivalent to Cartan's until Giaro, Kumpera, and Ruiz discovered the first singularity in 1978 (\cite{GKR}).

Jean (\cite{J}) studied the kinematic model of a car pulling $N$ trailers, a system which is locally universal for Goursat distributions of corank $N+1$.  He developed a geometric stratification given by regions in the configuration space of the model in terms of critical angles.  He also derived recurrence relations enabling one to compute the small growth vector of a Goursat germ for these geometric strata.  Montgomery and Zhitomirskii (\cite{MZ2}) introduced the Monster tower, a sequence of manifolds with distributions in which every Goursat germ occurs, as well as the Sandwich Lemma, allowing for Jean's strata to be recast in terms of positions of members of a canonical subflag of the Goursat flag.  Mormul (\cite{Mo1}) labelled the strata from \cite{MZ2} by words in the letters GST, which became the RVT code in \cite{MZ1}.  He solved Jean's relations in terms of the derived vector of the small growth vector, allowing for the calculation of the RVT code from the small growth vector.
  In \cite{MZ1}, Montgomery and Zhitomirskii showed that Goursat germs correspond to finite jets of Legendrian curve germs, and that the RVT coding corresponds to the classical invariant in the singularity theory of planar curves: the Puiseux characteristic (see Section \ref{Puiseux} below).  They gave a recursive algorithm for computing the Puiseux characteristic from the RVT code.

    In short, the present contribution provides the dashed arrow in the following diagram:
    
    \[ \begindc{\commdiag}[60]
    \obj(0,1)[SGV]{\{SGV\}}
    \obj(1,1)[RVT]{\{RVT\}}
    \obj(1,0)[PC]{\{PC\}}
    \mor{SGV}{RVT}{}[+1, 10]
    \mor{RVT}{PC}{}
    \mor{PC}{RVT}{}
    \mor{SGV}{PC}{}[+1,1]
    \enddc
    \]
       Here, the arrow \{RVT\} $\longleftrightarrow$ \{PC\} was given in \cite{MZ1},  and the arrow \{SGV\} $\longrightarrow$ \{RVT\} was given in \cite{Mo1}.  The arrow \{RVT\} $\longrightarrow$ \{SGV\} was given recursively in \cite{J}, and later explicitly in \cite{Mo1}.

It is worth noting that now, with both mappings \{RVT\} $\longrightarrow$ \{SGV\} and \{SGV\} $\dashrightarrow$ \{PC\} having explicit presentations, the recursive mapping Pc of \cite{MZ1} 
(see Section \ref{Pc}) has become explicit as well.

\section{Background}
\subsection{Goursat distributions and small growth vectors}

Let $M$ be a smooth manifold and let $D \subset TM$ be any smooth distribution (subbundle).  Let $[D, D] + D$ be called the $\textit{Lie square}$ of $D$.  Iterate this squaring to obtain a chain, which we write as
    \[ D^{s} \subseteq D^{s-1} \subseteq \cdots \subseteq D^{i} \subseteq D^{i-1} \subseteq \cdots \]
    where $D^{s}=D$ and $D^{i-1}$ is the Lie square of $D^{i}$.  Note that $D^{i}$ may not, in general, have constant rank, and thus fail to be a distribution on $M$.
    
     A rank 2 distribution \emph{D} on a manifold $M$ of real dimension $n=s+2\geq 4$ is called \emph{Goursat} if corank $D^{i} =i$ for $i=0, \dots, s$.
    In this case, one has a $\textit{Goursat flag} \ \mathcal{F}$:
    \[ D^{s} \subset D^{s-1} \subset \cdots \subset D^{1} \subset D^{0}=TM. \]
    Note that when $D$ is Goursat, each member $D^{i}$ of the flag is itself a distribution, and a hyperplane in $D^{i-1}$.

    Given a Goursat distribution $D$, one can alternatively form the sequence $D_{i}=[D, D_{i-1}]+D_{i-1}$, where $D_{1}=D$ and $i \geq 1$.  It is not hard to show that this sequence will also eventually terminate.  That is, there exists an $r$ such that $D_{r}=TM$.  Thus, Goursat distributions are \emph{completely nonholonomic}.  For any $p \in M$, the least $r$ such that $D_r(p) = T_pM$ 
    is called the \emph{degree of nonholonomy} of $D$ at $p$.  Note that for a general completely nonholonomic distribution, the degree of nonholonomy depends on the base point $p$, and so it happens for the Goursat distributions. 
        For each $p \in M$, we define the \emph{small growth vector at p} to be the integer valued vector
    \[sgv(p)= \big(\text{dim} D_{1}(p),\ \text{dim} D_{2}(p), \dots ,\ \text{dim} D_{r(p)}(p)=n\big).  \]
    In the following we shall only be concerned with germs of Goursat distributions.

    While the small growth vector is the traditional object of interest, for completely nonholonomic distributions it is equivalent to the \emph{derived vector}, which will be more convenient for us to work with.
    
    \begin{definition}[\cite{Mo1}]
    \label{derdef}
    The \emph{derived vector} of a completely nonholonomic distribution germ consists of the multiplicities of the entries in the small growth vector at the reference point.
    \end{definition}
    For a Goursat distribution, the dimensions of the sequence $D_{i}$ grow by at most one at a time, so the multiplicities (the entries in the derived vector) are nonzero.  By convention, we omit the last multiplicity 1 from the derived vector.  For example, if we are given a small growth vector (2, 3, 4, 4, 5), then the corresponding derived vector is (1, 1, 2).  Similarly, given a derived vector (1, 1, 1, 3, 3), the corresponding small growth vector is (2, 3, 4, 5, 5, 5, 6, 6, 6, 7).  While it is not obvious, for Goursat distributions the derived vector is always non-decreasing (see Section 2 of \cite{Mo1}).

\subsection{Construction of the RVT code}\label{RVT}
An RVT code is a word in the letters R,V,T satisfying one spelling rule: the letter T cannot follow the letter R.  An RVT code represents an equivalence class of Goursat germs.  The construction of a Goursat germ's RVT code was implicit in \cite{MZ2}, and made explicit in \cite{Mo1}, where the letters G,S,T were used instead of R,V,T.  Beginning with a Goursat germ, one forms the Goursat flag, which possesses a canonical integrable subflag called the \emph{characteristic foliation} or \emph{Cauchy characteristic}.  The geometric relationship between the members of the two flags can be characterized as Regular, Vertical, or Tangent (or, alternatively, Generic, Singular, or Tangent) and one encodes this information into a word called the RVT code.  See Section 1.2 of \cite{Mo1} for details.

In \cite{MZ1}, a parallel definition of the RVT codes for Goursat 
germs was proposed using a tower of manifolds called the ``Monster Tower.''  This tower is Goursat universal: every Goursat germ occurs somewhere within the Monster.  Each point in the Monster Tower is assigned an RVT code, and the code of a Goursat germ at a reference point $p$ is that of $p$ itself.

The tower is constructed through a series of {\it Cartan 
prolongations}.  Begin with the manifold $M^0 = \mathbb{R}^2$ and the distribution $\Delta^0=T\mathbb{R}^2$.  The first prolongation is the fiber bundle
\[ M^1 = \bigcup_{p\in \mathbb{R}^2}  \mathbb{P}\Delta^0_p,\]
whose elements have the form $(p,l)$, where $p$ is a point in $\mathbb{R}^2$ and $l$ is a line in the tangent space $T_p\mathbb{R}^2$.  The distribution on $M^1$ is given by
\[ \Delta^1_{(p,l)}= (d\pi^1_0)^{-1}(l) \]
where $\pi^1_0 \colon M^1 \to M^0$ is the bundle projection.

Iterating the prolongation procedure gives a sequence of manifolds
 \[ M^i = \bigcup_{p\in M^{i-1}}  \mathbb{P}\Delta_p^{i-1}.\]
  Every point in $M^i$ has the form $(p,l)$, where $p$ is a point in $M^{i-1}$ and $l$ is a line in the distribution plane $\Delta_p^{i-1}$. The dimension of $M^i$ is thus $i+2$.  The bundle projection map $\pi^i_{i-1} \colon M^i \to M^{i-1}$ has fibers diffeomorphic to $\mathbb{P}\Delta_p^{i-1} \cong \mathbb{RP}^1 \cong S^1$.  The distribution on $M^i$ is given by
  \[ \Delta^i_{(p,l)}= (d\pi^i_{i-1})^{-1}(l). \]
  One verifies that each distribution $\Delta^i$ is rank 2 and Goursat.  The \textit{Monster Tower} is thus the sequence of circle bundles
  \[ \dots \rightarrow M^i \rightarrow M^{i-1} \rightarrow \dots \rightarrow M^1 \rightarrow M^0=\mathbb R^2 \]
   equipped with a Goursat distribution at each level. 
   
  By composing the projection maps $\pi^k _{k-1},\  \pi^{k-1}_{k-2}, \ldots, \pi^{i+1}_{i}$ we obtain projections $\pi_i^k : M^k \to M^i$, $i < k$. The horizontal curves at level $i$  
naturally prolong (i.e., lift) to  horizontal curves at level $k$. 
However, the curves coinciding with fibers of $\pi_i^{i+1}$, 
$i \ge 1$, are special -- they project down to points and 
are not prolongations of curves from below. They are called 
\emph{vertical} and can themselves be prolonged to (first order) 
tangency curves, then prolonged again to (second order)
tangency curves, and so on. Vertical curves and their prolongations are called 
{\it critical}. 
Thus, at each level $i\geq 2$ there are vertical 
directions, and, in addition, at each level $i \geq 3$ there 
are tangency directions different from the vertical direction. At any level, 
all the remaining (non-critical) horizontal directions are called \emph{regular}. 

Recall that a point $p$ at level $i$ has the form $(q,l)$ where $q \in M^{i-1}$ and $l$ is a line in $\Delta^{i-1}$.  We call $p$ a \emph{regular, vertical,} or \emph{tangency} point if the direction of $l$ is regular, vertical, or tangency, respectively.  Points which are vertical or tangency are called \emph{critical}.  Therefore, at each level $i\geq 3$ there are regular and vertical 
points, and at each level $i \geq 4$ there 
are regular, vertical, and tangency points.

Now, the RVT code of a point $p$ at level $k \ge 3$ is 
a word $(\omega_3\omega_4\dots \omega_k)$ in the letters {R, V, T} 
satisfying 

$ \indent \omega_i = \begin{cases}
R &  \text{if}\ \pi_i^k(p) \ \text{is a regular point,}\\
{V} & \text{if}\ \pi_i^k(p) \ \text{is a vertical point,}\\
{T} & \text{if}\ \pi_i^k(p) \ \text{is a tangency point.}
\end{cases}$

It follows that the only spelling rule for the codes 
of points is the absence of the sequence `RT' in the codes; all other codes are realizable. Thus, the general form
of an RVT code of a point is $R^{k-2}$ (not dealt with in the present paper), or else
 \begin{equation}\label{code} 
 R^{r_{v+1}}VT^{t_v}R^{r_v} \cdots VT^{t_2}R^{r_2}VT^{t_1}R^{r_1}. 
 \end{equation}
Let us pause momentarily to explain the notation in (\ref{code}).  Here, $v$ is the number of letters V in the RVT code.  We have thus partitioned our code into $v+1$ pieces separated by the letters V.  We write that the \textit{last} letter V in the code is followed by $t_1$ many letters T, and then $r_1$ many letters R.  We continue for $1 \leq i \leq v$ letting $t_i$ denote the number of letters T following the $i$th V \textit{from the right}, and $r_i$ denote the number of letters R following those letters T.  Let $r_{v+1}$ denote the number of letters R preceding\footnote{Note that here, as in \cite{MZ1}, we allow $r_{v+1}\geq 0$.  In \cite{Mo1}, where the letters R,V,T are replaced by G,S,T, respectively, the code always begins with two letters G.  Thus, the GST codes in \cite{Mo1} are two letters longer than the RVT codes here and in \cite{MZ1}.} the first letter V.  Finally, superscripts in (\ref{code}) denote multiplicities.

The germ of a horizontal curve passing 
through a reference point $p$ is called {\it regular} when 
it is immersed and tangent neither to the vertical 
nor tangency direction at $p$. 
It is a central and deep fact in \cite{MZ1} that each germ of a non-constant 
well-parametrized analytic plane curve becomes regular after finitely many
prolongations and stays regular in 
 subsequent prolongations.  The least number of prolongations needed to regularize the curve is called the \emph{regularization level} $k$, and the $k$-fold prolongation of the original curve $\gamma$ is called the \emph{regularization} of $\gamma$.  This result bridges two seemingly 
distant areas: Goursat geometry and the singularity theory of  plane curves. 
In \cite{MZ1}, the first of two proofs of this deep fact is based on the Puiseux 
characteristic of a singular plane curve. 
It is worth noting that the reference point at level $k$, 
hit by the regularized curve, is still 
critical. Only its `son' at level $k + 1$ (hit by the 
prolongation of the newly regularized curve) and his sons will be regular.

Now the germ of a plane curve $\gamma$ with regularization level $k\geq 3$ is assigned an RVT code as well: it is the code of the (critical) 
reference point hit by the regularization of $\gamma$. 
Therefore, the codes of plane curves always end 
with V or T; 
such codes are called \emph{critical} (codes consisting solely of letters V and T are called \emph{entirely critical}). That is, 
in the notation from (\ref{code}), the equation $r_1 = 0$ holds for the codes of curves. 
The mapping Pc of \cite{MZ1}, one of the key players in the 
present work, acts precisely on the critical codes.

Lastly, note that if the original curve 
germ, or its first prolongation, is already immersed, 
then the curve's code is undefined. The further prolongations 
of such curves hit (and completely exhaust) the simplest 
points of the Monster: the so-called Cartan
points, or jet-like points (see Section \ref{intro}). At each level $k \ge 2$ 
these points populate the only open dense stratum $R^{k-2}$ (the entire $M^2$ when $k = 2$) which remains outside the field 
of interest of the present paper. \\ \\

We now recall the construction of the RVT code from the derived vector (see \cite{Mo1}). Suppose we are given a Goursat germ whose derived vector (see Definition \ref{derdef}) is
\[ der=(\underbrace{M_1, \ M_1, \dots, M_1}_{m_1}, \ \underbrace{M_2, \ M_2, \dots, M_2}_{m_2}, \dots, \underbrace{M_{v+1}, \ M_{v+1}, \dots, M_{v+1}}_{m_{v+1}}),  \]
with $M_1<M_2<\cdots<M_v<M_{v+1}$.

 Then $v$ turns out to be the number of letters V in the RVT code of this germ, which, therefore, has the form of (\ref{code}):
 \[ R^{r_{v+1}}VT^{t_v}R^{r_v} \cdots VT^{t_2}R^{r_2}VT^{t_1}R^{r_1}. \]
 Mormul derived the following relations for ascertaining the multiplicities $r_i$ and $t_i$.  See Theorem 3.5 in \cite{Mo1}. 
 One has: \\
    $\indent r_{v+1}=m_{v+1}-1 $ \\
    $\indent t_1=M_2 -2$ \\
    $\indent r_1=m_1-M_2 \geq 0$. \\
    For $2\leq i \leq v $ we have:
    
    \textbf{Case 1}:  $M_i$ divides $M_{i+1}$.  Then 
    \begin{align}
     \label{i} t_i&= \frac{M_{i+1}}{M_i} -2 \\
     \label{ii} r_i&= m_i-t_i -1>0. \\
    \intertext{\indent \textbf{Case 2}: $M_i$ does not divide $M_{i+1}$.  Then }
     \label{iii} t_i&= m_i-1 \\
     \label{iv} r_i&= 0.
     \end{align}

\subsection{The Puiseux characteristic}\label{Puiseux}
    Suppose $\gamma \colon (\mathbb{R}, 0) \to \mathbb{R}^2$ is the parametrization of an analytic plane curve germ.  We say $\gamma$ is  \textit{badly-parametrized} if there exist analytic germs $\mu \colon (\mathbb{R},0) \to \mathbb{R}^2$ and $\phi \colon (\mathbb{R},0)\to (\mathbb{R},0)$ such that $d\phi/ dt(0)=0$ and $\gamma=\mu \circ \phi$.  Otherwise, $\gamma$ is called \textit{well-parametrized}.  If $\gamma$ is well-parameterized and not immersed then we may define its Puiseux characteristic, an invariant with respect to the RL-equivalence of curve germs.  Up to RL-equivalence, $\gamma$ has the form
\[ \gamma(t)=(t^m, \sum_{k \geq m} a_k t^k) \]
where $m \geq 2$.

The definition of the Puiseux characteristic is the following.
Let $\lambda_0=e_0=m$.  Then define inductively for $j \geq 0$
\[ \lambda_{j+1} = \text{min} \{ k \ | \ a_k \neq 0,\ e_j \nmid k\}, \quad \quad e_{j+1}= \text{gcd}(e_j, \lambda_{j+1}) \]
until we first obtain a $g$ with $e_g=1$.  Then the vector $[\lambda_0 ; \lambda_1, \dots, \lambda_g]$ is called the \textit{Puiseux characteristic} of $\gamma$.

The Puiseux characteristic is the fundamental invariant in the singularity theory of plane curves.  In \cite{W}, Proposition 4.3.8 shows that it is equivalent to at least seven other classical invariants.

Here, as in \cite{MZ1}, we restrict our attention to Puiseux characteristics satisfying
\begin{equation}\label{lambdas}
\lambda_1 > 2\lambda_0.
\end{equation}
This is a normalization condition, and no equivalence classes of Legendrian curves are excluded by its imposition.  The Puiseux characteristic is an invariant with respect to RL-equivalence, but not a complete invariant.  For example, $(t^2, t^5)$ and $(t^3, t^5)$ have different Puiseux characteristics, but their first prolongations are equivalent Legendrian curves.  The restriction (\ref{lambdas}) resolves this ambiguity.

\subsection{The map Pc}\label{Pc}
Here we recall the definition of the map Pc constructed in Section 3.8.4 of \cite{MZ1}.  Given a critical RVT code $(\alpha)$, this map yields a Puiseux characteristic $\text{Pc}(\alpha)$ satisfying (\ref{lambdas}).  The relationship between $(\alpha)$ and $\text{Pc}(\alpha)$ is given in Theorem 3.23 of \cite{MZ1}.  The map is constructed recursively as follows.

First define the two maps
\begin{align*}
\mathbb E_T &\colon (n_1, n_2) \mapsto (n_1, n_1+n_2) \\
\mathbb E_V &\colon (n_1, n_2) \mapsto (n_2, n_1+n_2).
\end{align*}
Then, for an entirely critical code
\[ (\omega)=(\omega_1, \dots, \omega_m), \quad \omega_i \in\{V,T\}
\]
define $\mathbb E_{\omega}$ to be the composition
\[\mathbb E_{\omega}=\mathbb E_{\omega_1} \circ \mathbb E_{\omega_2} \circ \cdots \circ \mathbb E_{\omega_m}.
\]
Next, we note that any critical RVT code $(\alpha)$ has one of the two following forms: 

A. $(\alpha)=(R^s\omega)$, where $s\geq 0$ and $(\omega)$ is an entirely critical RVT code;

B. $(\alpha)=(\beta R^s\omega)$, where $s\geq 1$, $(\beta)$ is a critical RVT code, and $(\omega)$ is an entirely critical RVT code. 

In case A, let $(a,b)=\mathbb E_{\omega}(1,2)$.  Then
\[ \text{Pc}(\alpha)= [\lambda_0;\lambda_1], \quad \lambda_0=a,\ \lambda_1=sa+a+b.
\]

In case B, let $(a,b)=\mathbb E_{\omega}(1,2)$ and $\text{Pc}(\beta)= [\tilde \lambda_0;\tilde \lambda_1, \dots, \tilde \lambda_{g-1}]$.  Then
\[ \text{Pc}(\alpha)= [\lambda_0;\lambda_1, \dots, \lambda_g], \]
where
\begin{align*}
\lambda_i&=a\tilde \lambda_i \quad \text{for} \quad 0 \leq i \leq g-1\\
\lambda_g &= a(\tilde \lambda_{g-1}+s-1)+b-a.
\end{align*}

\section{Main Result}
As explained in Section 1, here we compose the formulas presented in \cite{Mo1} and the algorithm in \cite{MZ1}, yielding a formula for the Puiseux characteristic of the plane curve corresponding to a Goursat germ with given small growth vector.  This formula turns out to be simpler than either of the two from which it was derived, suggesting a deeper geometric link between singularities of plane curves and singular Goursat distributions.  The problem solved herein was first proposed in \cite{MZ1} as Question 9.19, part 3, and was asked again in \cite{Mo2} in the Afterword.

\subsection{Main Theorem}
Suppose we are given a Goursat germ whose derived vector (see Definition \ref{derdef}) is
\[ der=(\underbrace{M_1, \ M_1, \dots, M_1}_{m_1}, \ \underbrace{M_2, \ M_2, \dots, M_2}_{m_2}, \dots, \underbrace{M_{v+1}, \ M_{v+1}, \dots, M_{v+1}}_{m_{v+1}}),  \]
with $M_1<M_2<\cdots<M_v<M_{v+1}$.
Consider the set $S=\{ M_i | \ M_{i-1} \ \text{divides}\ M_i \}$.  Let $g=|S|$.  For $1\leq j\leq g$, let $N_1, N_2, \dots, N_g$ denote the elements of $S$ in decreasing order.  We always have $N_g=M_2$, since $M_1=1$. For $1\leq j\leq g$ let $M_{k_j}=N_j$.  

\begin{theorem}
\label{maintheorem}
The corresponding Puiseux characteristic is
\[ [\lambda_0; \lambda_1, \dots, \lambda_{g}] \]
where
\begin{align}
\label{eqn1}
\lambda_0&= M_{v+1} \\
\label{eqn2}
\lambda_j &=\sum_{i\geq k_j} m_i M_i+M_{k_j} +M_{k_j-1}
\end{align}
for $1\leq j\leq g$.
\end{theorem}

\subsection{Example}
Suppose $der=(1, 1, 2, 2, 2, 2, 2, 2, 4, 6, 6, 6, 18, 24, 24)$.  The associated RVT code is $RVVTRVVRRRRRV$.  Note that $\lambda_0=M_{v+1}=M_{6}=24$.  We also have
$S=\{18, 4, 2\}$, and therefore $g=3$.  Then write $S=\{18, 4, 2 \}=\{N_{1}, N_2, N_3 \}=\{M_5, M_3, M_2 \}$ so that $k_1=5,\ k_2=3,$ and $k_3=2$.
Finally, we compute
\begin{align*}
\lambda_1&= \sum_{i\geq 5} m_iM_i+M_5+M_4=90 \\
\lambda_2&= \sum_{i\geq 3} m_iM_i+M_3+M_2=94 \\
\lambda_3&= \sum_{i\geq 2} m_iM_i+M_2 +M_1=103. \\
\end{align*}
The Puiseux characteristic is thus
\[ [24; 90, 94, 103]. \]

\subsection{Example}
This example is very similar to the previous, and the subtle differences should provide room for comparison.  It appears in \cite{MZ1} as Example 3.28.

Suppose $der=(1, 1, 2, 2, 2, 2, 2, 2, 4, 6, 6, 6, 6, 18, 24, 24)$.  The associated RVT code is $RVVTRRVVRRRRRV$.  Note that we have the same values of $g, \ N_j,\ k_j,$ and $M_1, M_2, \dots, M_{v+1}$  as in the previous example.
Thus, we compute
\begin{align*}
\lambda_1&= \sum_{i\geq 5} m_iM_i+M_5+M_4=90 \\
\lambda_2&= \sum_{i\geq 3} m_iM_i+M_3+M_2=100 \\
\lambda_3&= \sum_{i\geq 2} m_iM_i+M_2 +M_1=109. \\
\end{align*}
The Puiseux characteristic is thus
\[ [24; 90, 100, 109]. \]

\subsection{Remark}
We are implicitly assuming that the underlying RVT class is critical (ends with V or T).  Then the associated planar curve germ is necessarily non-immersed.  We only discuss the Puiseux characteristic for these singular planar curve germs, since any immersed planar curve germ has normal form $(t, 0)$.  This restriction agrees with the domain of the map Pc given in Section \ref{Pc}.  In terms of the derived vector, according to Section \ref{RVT}, we must assume that $m_1=M_2$.

\section{Proof of Theorem}
The theorem is proved by induction on $g$.  For readability, we break the proof into three subsections.  In the first, we prove a useful lemma.  In the second, we verify the base case $g=1$.  In the third, we complete the inductive step.

\subsection{Lemma}
 The following lemma makes use of the basic bricks $A_i$ constructed in Section 3.1 of \cite{Mo1}.  These are integers from which the entries $M_i$ in the derived vector are built, and the two actually coincide in some cases -- for the precise relationship, see Theorems 3.3 and 3.4 in \cite{Mo1}.  The $A_i$ (and subsequently the $M_i$) depend only on the parameters $t_1, t_2, \dots, t_v$; the multiplicity $m_i$ depends on both $t_i$ and $r_i$.  The bricks are constructed as follows:
\begin{align*}
A_1&=1\\
A_2&= 2+t_1 \\
A_i&= A_{i-2}+A_{i-1}(1+t_{i-1}) \quad  \text{for} \quad 3\leq i \leq v+1.
\end{align*}

\begin{lemma}\label{cor}
Let $(\omega)=VT^{t_N}\cdots VT^{t_2}VT^{t_1}$ for $N \geq 2$. Let $(a,b) = \mathbb E_{\omega}(1,2)$.  Then
\begin{align*}
 a&=A_{N+1}  \\
 b&= A_1+A_2+\sum^N_{i=2} (1+t_i)A_i.
\end{align*}

\end{lemma}

\begin{proof}
The lemma follows from the following two observations:
\begin{align}
\label{observ1}A_1+A_2+\sum^N_{i=2} (1+t_i)A_i= A_N+A_{N+1} \\
\intertext{and}
\label{observ2}\mathbb E_V \mathbb E_{T}^{t_N} \cdots \mathbb E_V \mathbb E_{T}^{t_2} \mathbb E_V \mathbb E_{T}^{t_1}(1,2)=(A_{N+1}, A_N+A_{N+1}).
\end{align}
Both observations can be easily verified for $N=2$.  Assuming Equation (\ref{observ1}) holds for $N$, we find
\begin{align*}
A_1+A_2+\sum^{N+1}_{i=2} (1+t_i)A_i&= A_1+A_2+\sum^N_{i=2} (1+t_i)A_i +(1+t_{N+1})A_{N+1}\\
&= A_N+A_{N+1}+(1+t_{N+1})A_{N+1} \\
&=A_{N+1}+A_{N+2},
\end{align*}
so (\ref{observ1}) is proved by induction.  Similarly, assuming Equation (\ref{observ2}) holds for $N$, we find
\begin{align*}
\mathbb E_V \mathbb E_{T}^{t_{N+1}}\mathbb E_V \mathbb E_{T}^{t_N} \cdots \mathbb E_V \mathbb E_{T}^{t_2} \mathbb E_V \mathbb E_{T}^{t_1}(1,2) &= 
\mathbb E_V \mathbb E_{T}^{t_{N+1}}(A_{N+1}, A_N+A_{N+1}) \\
&= \mathbb E_V (A_{N+1}, t_{N+1}A_{N+1}+A_N+A_{N+1}) \\
&=\mathbb E_V(A_{N+1}, A_{N+2}) \\
&= (A_{N+2}, A_{N+1}+A_{N+2}),
\end{align*}
so (\ref{observ2}) is proved by induction.  This completes the proof of the lemma.
\end{proof}

\subsection{Base case of the induction}
For the base case, assume $g=1$. Note that $k_1=k_g=2$.  Now, according to Section \ref{RVT} we have $r_i=0$ for $i=2, \dots, v$, so the associated RVT code is of the form $R^{s}\omega$, with $(\omega)$ entirely critical (containing no letters R).  We also have $m_{v+1}=s+1$ and $m_i=1+t_i$ for $i=2, \dots, v$.  Then Lemma \ref{cor}, with $N=v$, gives
\begin{align*}
a&=A_{v+1} \\
b&= A_1+A_2+\sum^v_{i=2} (1+t_i)A_i \\
&=A_1+A_2+\sum^v_{i=2} m_iA_i.
\end{align*}
Then, according to Section \ref{Pc}, we have that the Puiseux characteristic is $ [\lambda_0; \lambda_1]$, with
\begin{align*}
\lambda_0 &= a =A_{v+1} \\
\lambda_1 &= sa+a+b \\
&=(m_{v+1}-1)A_{v+1} + A_{v+1} +A_1+A_2+\sum^v_{i=2} m_iA_i \\
&=A_1+A_2+\sum^{v+1}_{i=2} m_iA_i.
\end{align*}
However, according to Theorem 3.3 from \cite{Mo1}, in the case $g=1$ we have that $A_i=M_i$ for all $i$.  Thus, the Puiseux characteristic is $ [\lambda_0; \lambda_1]$ with
\begin{align*} 
\lambda_0 &=M_{v+1} \\
 \lambda_1 &=M_1+M_2+\sum^{v+1}_{i=2} m_iM_i
\end{align*}
as desired.

\subsection{Inductive step}
We now assume Theorem \ref{maintheorem} for derived vectors satisfying $|S|<g$, and prove that the result must hold for derived vectors with $|S|=g$.  We may assume $g>1$.

Begin with an arbitrary derived vector 
\[ der=(\underbrace{M_1, \ M_1, \dots, M_1}_{m_1}, \ \underbrace{M_2, \ M_2, \dots, M_2}_{m_2}, \dots, \underbrace{M_{v+1}, \ M_{v+1}, \dots, M_{v+1}}_{m_{v+1}})  \]
for which $|S|=g$.
The idea is to truncate the associated RVT code $(\alpha)$ after the last occurring letter R.  Our inductive hypothesis will then apply to the derived vector associated to the truncated code, and we can reconstruct the Puiseux characteristic of the original derived vector from here.  

To this end, we must give special attention to the entry $N_{g-1}=M_{k_{g-1}}$ in \emph{der}.  For notational purposes, we set $q=k_{g-1}$.  Then by assumption we have that $M_{q-1}$ divides $M_q$, and $M_{q-1}$ is the smallest such entry (besides $M_1=1$).  

The relations in Section \ref{RVT} imply that we can write
    \[(\alpha)=(\beta R^s \omega)
    \]
    with
        \begin{align*}
        s&>0 \\
        \omega &= \text{enitrely critical RVT string}\\
        \beta &= \text{critical RVT code}.
        \end{align*}
More explicitly, we write
\[(\alpha)=\underbrace{R^{r_{v+1}}VT^{t_v}R^{r_v} \cdots VT^{t_{q-1}}}_{\beta}R^{r_{q-1}}
\underbrace{VT^{t_{q-2}}VT^{t_{q-3}}\cdots VT^{t_{1}}}_{\omega}
\]
and note that $s= r_{q-1}=m_{q-1}-\frac{M_{q}}{M_{q-1}}+1$, by Equations (\ref{i}) and (\ref{ii}).

Now, $\beta$ is our truncated code and we will adorn all data concerning ($\beta$) with tildes to distinguish them from those of ($\alpha$).  In particular, we write $[\tilde{\lambda}_0; \tilde{\lambda}_1, \dots, \tilde{\lambda}_{g-1}]$ for the Puiseux characteristic of $(\beta)$ and
    \[(\tilde{M}_1^{\tilde{m}_1}, \tilde{M}_2^{\tilde{m}_2}, \dots, \tilde{M}_{\tilde{v}+1}^{\tilde{m}_{\tilde{v}+1}})\]
    for the derived vector.  Note that $\tilde{g}$ is indeed equal to $g-1$ by virtue of the recursive description of the mapping Pc in \cite{MZ1} (see Section \ref{Pc}).  

With this setup, we begin the calculations.  Applying Lemma \ref{cor} with $N=q-2$ we find that $(a,b)=\mathbb E_{\omega}$ is given by 
\begin{align*}
a&=A_{q-1} \\
b&= A_1+A_2+\sum^{q-2}_{i=2} (1+t_i)A_i.
\end{align*}
But according to Theorem 3.4 in \cite{Mo1}, we have\footnote{In \cite{Mo1}, the quantity $n_1$ tells us the position of the last occurring letter R in a critical RVT code.  Here, we have $n_1=q-2$, which means $r_{q-1}\neq 0$ and $r_i=0$ for $i<q-1$; in other words, $\omega$ is the largest entirely critical string at the tail of $(\alpha).$} $A_i=M_i$ for $1\leq i \leq q-1$.  Also, by Equations (\ref{ii})-(\ref{iv}), we have $m_i=1+t_i+r_i$ for all $i\geq 2$.  But $r_1=r_2=\cdots=r_{q-2}=0$, so for $2\leq i \leq q-2$ we have $1+t_i=m_i$.  Thus,
\begin{align*}
a&=M_{q-1} \\
b&= M_1+M_2+\sum^{q-2}_{i=2} m_iM_i.
\end{align*}

Next, we need to relate the derived vector of $(\alpha)$ to that of $(\beta)$.  Proposition 1 from \cite{Mo2} gives the following relations:
\begin{align*}
\tilde m_1 &= \frac{M_q}{M_{q-1}} \\
\tilde m_i &= m_{q-2+i} \quad \text{for} \quad 2\leq i \leq \tilde v+1 \\
\tilde M_i &= \frac{M_{q-2+i}}{M_{q-1}} \quad \text{for} \quad 1\leq i \leq \tilde v+1  \\
\tilde v &= v-q+2.
\end{align*}
Also, according to Section \ref{Pc}, we know how to compute the Puiseux characteristic of $(\alpha)$ from that of $(\beta)$:
\begin{align*}
\lambda_j&=a\tilde \lambda_j \quad \text{for} \quad 0\leq j \leq g-1 \\
\lambda_g &= a(\tilde \lambda_{g-1}+s-1)+b-a.
\end{align*}
But by induction we may assume
    \begin{align*}
    \tilde{\lambda}_0&= \tilde{M}_{\tilde{v}+1} \\
    \tilde{\lambda}_j &=\sum_{i\geq \tilde{k}_j} \tilde{m}_i \tilde{M}_i+\tilde{M}_{\tilde{k}_j} +\tilde{M}_{\tilde{k}_j-1}
    \quad \text{for} \quad 1\leq j\leq g-1.
    \end{align*}

Putting this all together, we now compute the desired Puiseux characteristic of $der=der(\alpha)$ in three stages.  First we compute $\lambda_0$, then $\lambda_j$ for $1\leq j \leq g-1$, then finally $\lambda_g$.

First, we easily find 
\begin{align*}
\lambda_0&=a\tilde \lambda_0 \\
&= M_{q-1}\frac{M_{q+\tilde v-1}}{M_{q-1}}\\
&= M_{v+1}.
\end{align*}
Next, for $1\leq j \leq g-1$ we have
\begin{align*}
\lambda_j &= a \tilde \lambda_j \\
&= M_{q-1}\Big(\sum_{i\geq \tilde{k}_j} \tilde{m}_i \tilde{M}_i+\tilde{M}_{\tilde{k}_j} +\tilde{M}_{\tilde{k}_j-1}\Big)\\
&= \sum_{i\geq \tilde{k}_j} m_{q-2+i}M_{q-2+i}+ M_{q+\tilde k_j-2} + M_{q+\tilde k_j-3} \\
&= \sum_{i\geq q+\tilde{k}_j-2}m_{i}M_{i}+ M_{q+\tilde{k}_j-2} + M_{q+\tilde{k}_j-3} \\
&=\sum_{i\geq k_j} m_i M_i+M_{k_j} +M_{k_j-1}.
\end{align*}
The last equality comes from the general fact that for $1 \leq j \leq g-1$ we have
\begin{equation}\label{index}
q+\tilde k_j-2=k_j.
\end{equation}
To see this, recall that $k_j$ is defined so that $N_j=M_{k_j}$ is the $j$th smallest entry in $der$ which is divisible by the preceding entry.  Since this observation applies to the derived vectors of both $(\alpha)$ and $(\beta)$, we find that
    \begin{align*}
    \tilde{M}_{\tilde{k}_j-1}\ \text{divides}\ \tilde{M}_{\tilde{k}_j} &\Leftrightarrow  \frac{M_{q+\tilde{k}_j-3}}{M_{q-1}}\ \text{divides} \ \frac{M_{q+\tilde{k}_j-2}}{M_{q-1}} \\
    &\Rightarrow M_{q+\tilde{k}_j-3}\ \text{divides}\  M_{q+\tilde{k}_j-2} \\
    &\Rightarrow M_{q+\tilde{k}_j-2} =M_{k_j} \\
    &\Rightarrow q+\tilde{k}_j-2 =k_j.
    \end{align*}

Lastly, we compute $\lambda_g$.  From above, we know 
\begin{align*}
a\tilde \lambda_{g-1}=\lambda_{g-1}&=\sum_{i\geq k_{g-1}} m_i M_i+M_{k_{g-1}} +M_{k_{g-1}-1} \\
&=\sum_{i\geq q} m_i M_i+M_{q} +M_{q-1}.
\end{align*}
Whence
\begin{align*}
\lambda_g&= a(\tilde \lambda_{g-1}+s-1)+b-a \\
&=\Big(\sum_{i\geq q} m_i M_i+M_{q} +M_{q-1} +m_{q-1}M_{q-1}-M_q\Big)+ \Big(M_1+M_2+\sum^{q-2}_{i=2} m_iM_i\Big)-M_{q-1} \\
&=M_1+M_2+\sum_{i\geq 2} m_iM_i.
\end{align*}
Since we always have $k_g=2$, this is the desired result.

\begin{flushright}
    $\blacksquare$
    \end{flushright}

\subsection*{Acknowledgments}
The author would like to thank Richard Montgomery, Wyatt Howard, and Alex Castro for helpful discussion, advice, and most importantly, motivation.  The author is also grateful to the referee, whose suggestions greatly improved the overall presentation of this paper.

\end{document}